\documentclass{birkmult}

\newtheorem{thm}{Theorem}[section]
 \newtheorem{cor}[thm]{Corollary}
 
 \newtheorem{prop}[thm]{Proposition}
 \theoremstyle{definition}
 
 \theoremstyle{remark}
 \newtheorem{rem}[thm]{Remark}
 \newtheorem*{ex}{Example}
 \numberwithin{equation}{section}

\usepackage{dsfont}
\usepackage{tikz}
\usepackage{amssymb}
\usepackage{amsmath}
\usepackage{enumerate}
\usepackage[all]{xy}
\usepackage{youngtab}
\usepackage{mathtools}
\usepackage{framed}
\usepackage{mathabx,epsfig}

\DeclareMathOperator\C{\mathbb C}
\DeclareMathOperator\Z{\mathbb Z}
\DeclareMathOperator\N{\mathbb N}

\DeclareMathOperator\T{\mathbb T}
\DeclareMathOperator\csm{c^{sm}}
\DeclareMathOperator\ssm{s^{sm}}
\DeclareMathOperator\mC{mC}
\DeclareMathOperator\mS{mS}

\DeclareMathOperator\pt{pt}

\DeclareMathOperator\Fl{\mathcal F\!{}_\lambda}

\DeclareMathOperator\h{\hbar}
\DeclareMathOperator\E{\mathbb E}
\DeclareMathOperator\HH{\mathbb H}
\DeclareMathOperator\K{\mathbb K}
\DeclareMathOperator\Sym{Sym}
\DeclareMathOperator\tb{\boldsymbol t}
\DeclareMathOperator\mub{\boldsymbol \mu}
\DeclareMathOperator\zb{\boldsymbol z}
\DeclareMathOperator\Il{{\mathcal I}_\lambda}
\DeclareMathOperator\id{id}
\DeclareMathOperator\LR{c}
\DeclareMathOperator\Ell{Ell}
\DeclareMathOperator\Gr{Gr}
\DeclareMathOperator\TGr{T^*\!\Gr}
\DeclareMathOperator\spa{span}
\DeclareMathOperator\PPP{\mathbb P}
\DeclareMathOperator\One{\mathds 1}

\newcommand*{\defeq}{\stackrel{\text{def}}{=}}

\title[$\h$-deformed Schubert calculus]{$\h$-deformed Schubert calculus in equivariant \\cohomology, K-theory, and elliptic cohomology}

\author{Rich\'ard Rim\'anyi}
\address{Department of Mathematics, University of North Carolina at Chapel Hill, USA}
\email{rimanyi@email.unc.edu}



\thanks{The author is grateful to A. Weber for several helpful comments on this paper, and to L.~M.~Feh\'er and A. Varchenko for many useful conversations on this topic. The author is supported by a Simons foundation grant. }

\begin{document}

\begin{abstract}
In this survey paper we review recent advances in the calculus of Chern-Schwartz-Mac\-Pher\-son, motivic Chern, and elliptic classes of classical Schubert varieties. These three theories are one-parameter ($\h$) deformations of the notion of fundamental class in their respective extraordinary cohomology theories. Examining these three classes in conjunction is justified by their relation to Okounkov's stable envelope notion. We review formulas for the $\h$-deformed classes originating from Tarasov-Varchenko weight functions, as well as their orthogonality relations. As a consequence, explicit formulas are obtained for the Littlewood-Richardson type structure constants.    
\end{abstract}

\maketitle


\ \hfill Dedicated to Andr\'as N\'emethi, on the occasion of his 60th birthday.

\section{Introduction}

A basic structure of traditional Schubert calculus is the cohomology ring of a homogeneous space $X$, together with a distinguished basis. The elements of the distinguished basis are associated with the geometric subvarieties (called Schubert varieties) of $X$. The first objects to study are the structure constants of the ring with respect to the distinguished basis. These structure constants satisfy various positivity, stabilization, saturation, and other properties, and can be related with other mathematical fields, such as combinatorics, representation theory, integrable systems.

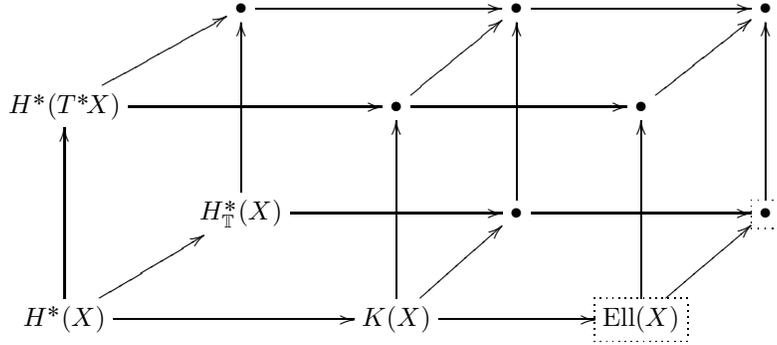
\begin{figure}\label{eq:kocka}
\xymatrix{
&  \bullet\ar@{>}[rr] & &  \bullet\ar@{>}[rr] & & \bullet \\
H^*(T^*\!X)\ar@{>}[ru]\ar@{>}[rr] & &  \bullet\ar@{>}[ru] \ar@{>}[rr]& &  \bullet\ar@{>}[ru] & \\
&  H^*_{\T}(X) \ar@{>}[rr] \ar@{>}[uu]& & \bullet \ar@{>}[rr]\ar@{>}[uu]& & \bullet\ar@{>}[uu] \\
H^*(X) \ar@{>}[uu] \ar@{>}[rr] \ar@{>}[ru]& &  K(X)\ar@{>}[rr]\ar@{>}[ru]\ar@{>}[uu]  & & \Ell(X)\ar@{>}[ru]\ar@{>}[uu] &  
\save "4,5"."4,5"*[F.]\frm{} \restore
\save "3,6"."3,6"*[F.]\frm{} \restore
}
\caption{Three ``orthogonal'' directions to generalize classical Schubert calculus}
\end{figure}

In this paper we survey some generalizations of this traditional setup, organized as vertices of the diagram in Figure \ref{eq:kocka}.   
The traditional setup, mentioned in the first paragraph, is the bottom left corner, the cohomology ring of the homogeneous space $X$. Going from the front face to the back face of the diagram represents the change from ordinary ``cohomologies'' to equivariant ones. Equivariant cohomologies of $X$ take into account the geometry of $X$ together with the natural group (say torus) action on it. Due to some techniques that only exist in equivariant theories it can be an {\em easier} theory to work with. Nevertheless, formulas in non-equivariant theories can be recovered from equivariant ones by plugging in $0$ (or $1$, depending on conventions) in the equivariant variables.  

Stepping one step to the right on the diagram from cohomology we arrive at K theory. K theoretic Schubert calculus, ordinary or equivariant, has been studied extensively, see \cite{LSS} and references therein. Stepping one further step to the right we get to elliptic cohomology, ordinary or equivariant: the two vertices in the diagram that are in dotted frames. We framed these vertices in the diagram because these settings lack the notion of a well-defined distinguished basis, which was present in $H^*$, $H^*_{\T}$, $K$, and $K_{\T}$. Namely, it turns out that in elliptic cohomology the notion of fundamental class depends on choices \cite{BE}. There are important results in these settings (e.g. \cite{GR, LZ, LZZ} and references therein) that follow from making certain choices (of a resolution, or a basis in a Hecke algebra).   

\begin{rem} There are other extraordinary cohomology theories, for example the (universal) complex cobordism theory; their position would be further to the right on the diagram. Yet, we restrict our attention to the ones depicted in Figure 1, as their associated formal group law is an algebraic group. 
\end{rem}

The focus of this paper is the rest of the diagram, namely the top face. Besides pioneering works, e.g. \cite{PP, AM1, AM2}, this direction of generalization is very recent. 

There are two ways of introducing this direction of generalization. One way is that we study Schubert calculus not on the homogeneous space $X$ but on its cotangent bundle $T^*\!X$, using its extra holomorphic symplectic structure. Although $X$ and $T^*X$ are homotopically equivalent so their cohomology rings are isomorphic, their geometric subvarieties are different (in $T^*X$ one studies the conical Lagrangian cycles).  We will not touch upon this interpretation, except in notation: In the diagram above we indicate this generalization by replacing $X$ with $T^*X$.

The other way of explaining the step from the bottom face to the top face is that we change the {notion of the characteristic class} that we associate with Schubert cells. While traditionally (on the bottom face) we associate the {\em fundamental class} to the Schubert variety, in this generalization we associate a particular one-parameter deformation of the fundamental class. The parameter will be denoted by $\h$, and the class will be called the $\h$-deformed Schubert class. In cohomology this $\h$-deformed class was conjectured/discovered by Grothendieck, Deligne, MacPherson \cite{M}, and is called the {\em Chern-Schwartz-MacPherson class}.\footnote{in the classical CSM literature the parameter $\h$ is not indicated, because it can be recovered from the grading in $H^*$.} Its equivariant theory is worked out by Ohmoto \cite{W1, O1, O2}. In K theory the $\h$-deformed class
\footnote{in most of the literature the letter $y$ is used for $\h$ in K theory, to match the classical notion of $\chi_y$-genus.} 
was defined by Brasselet-Sch\"urmann-Yokura \cite{BSY} under the name of {\em motivic Chern class}. The equivariant version is defined in \cite{FRW1,AMSS2}.\footnote{The CSM, motivic Chern, and elliptic classes were {\em not} discovered as $\h$-deformations of the notion of the fundamental class, but as generalizations of the notion of total Chern class for singular varieties with covariant functoriality; their interpretation as $\h$-deformations of the fundamental class suggested in this paper is post-factum. Thus, the present paper is a re-interpretation of the story of characteristic classes of singular varieties from the mid-70s to the present.}

The most recent discovery is the definition of the (ordinary or equivariant) $\h$-deformed elliptic class associated with a Schubert cell \cite{RW, KRW}, that is, Schubert calculus in the rightmost two vertices on the top face of the diagram. A pleasant surprise of such Schubert calculus is that the $\h$-deformed elliptic class does not depend on choices---the corresponding vertices in the diagram are not framed. While it is a general fact that one can recover the non-$\h$-deformed theory from the $\h$-deformed theory by plugging in an obvious value ($0,1,\infty$, depending on setup) for $\h$, it turns out that at such specialization the $\h$-deformed elliptic class has a singularity. This fact is another incarnation of the phenomenon mentioned above that the non-$\h$-deformed elliptic Schubert calculus depends on choices.

\begin{rem} \label{rem:Ok}
Let us comment on a principle that unifies the three $\h$-deformations in $H^*$, $K$, $\Ell$, which is actually the reason for the attention $\h$-deformed Schubert calculus is getting recently. In works of Okounkov and his co-authors Maulik, Aganagic \cite{MO, O, AO} (see also \cite{RTV1, GRTV, RTV2, RTV3, RTV4}), a remarkable bridge is built between quantum integrable systems and geometry. Via this bridge the (extraordinary) cohomology of a geometrically relevant space is identified with the Bethe algebra of a quantum integrable system.
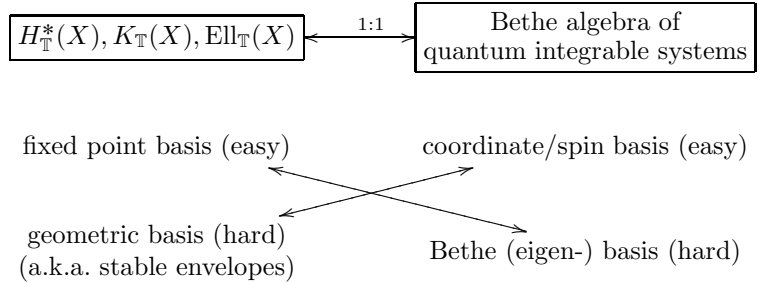
\[
\xymatrix @=1.5pc{
*+[F]{ H_{\T}^*(X),K_{\T}(X), \Ell_{\T} (X)} \ar@{<->}[rr]^{1:1}&&
 *+[F] \txt{Bethe algebra of \\ quantum integrable systems} \\
\txt{fixed point basis (easy)} \ar@{<->}[rrd]
&& \txt{coordinate/spin basis (easy)} \\
\txt{geometric basis (hard)\\ (a.k.a. stable envelopes)}  \ar@{<->}[rru]
 && \txt{Bethe (eigen-) basis (hard)}
}
\]
On both sides of the identification we have an ``easy'' and a ``hard'' basis, and the identification matches the easy basis of one side with the hard basis of the other side. The geometric basis that matches the spin basis of the Bethe algebra side is named the cohomological, K theoretic, and elliptic {\em stable envelope}. It is now proved that in type A Schubert calculus settings the three flavors of stable envelopes coincide (through some identifications, and convention matching) with the three $\h$-deformed Schubert classes: the CSM class, the motivic Chern class, and the elliptic class. This relation with quantum integrable system is the reason we denote the deformation parameter by $\h$.
\end{rem}

\begin{rem}
A fourth direction to generalize classical Schubert calculus is {\em quantum Schubert calculus}. While quantum cohomology and K theory (and possibly elliptic cohomology) are related with their $\h$-deformations see e.g. \cite{MO}, \cite[App.3]{RTV3}, we will not study them in this paper.
\end{rem}

\smallskip

The topic of $\h$-deformed Schubert calculus is rather fresh and the available literature is rather technical. Moreover, as explained in Remark \ref{rem:Ok} above, some of the existing literature is hidden in quantum integrable system papers. The goal of the present survey is twofold. On the one hand we want to give an accessible, well-motivated, and technicality-free presentation of the main achievement, what we call Main Theorem, see Section \ref{sec:main}. On the other hand we give precise formulas of the key ingredients (weight functions, inner products, orthogonality statements) consistent with usual Schubert calculus usage. While these formulas exist in some conventions in the literature, the conventions used there are optimized for some other purposes. Also, we tried to separate the complicated formulas (they are exiled to the penultimate section) from the main part of the paper where the idea is presented.

\section{Ordinary and equivariant cohomological Schubert calculus}
\subsection{Schubert classes and structure coefficients} \label{sec:SchLR}
Consider the compact smooth variety $\Gr(m,n)$, the Grassmannian of $m$-planes in $\C^n$. For an $m$-element subset $I$ of $[n]\defeq\{1,\ldots,n\}$ one defines the Schubert cell
\[
\Omega_I=\{V^m\subset \C^n : \dim(V^m \cap \C^q)=|\{i\in I:i\leq q \}|\ \forall q\},
\]
where $\C^1\subset \C^2\subset \ldots \subset \C^{n-1}\subset \C^n$ is the standard full flag. The collection of cohomological fundamental classes $[\overline{\Omega}_I]$ forms a basis in the cohomology ring of the Grassmannian, hence via 
\[
[\overline{\Omega}_I] \cdot [\overline{\Omega}_J] = \sum_K \LR_{I,J}^K \cdot [\overline{\Omega}_K]
\]
the structure coefficients (a.k.a. Littlewood Richardson coefficients) $\LR_{I,J}^K\in \Z$ are defined. 

To name an example, let us ``encode'' the subset $I=\{i_1<\ldots<i_m\}$ with the partition $(\lambda_1\geq \lambda_2 \geq \ldots\geq \lambda_m)$, by $\lambda_j=n-m-(i_j-j)$. With a slight abuse of notation let the Young diagram of $\lambda$ mean the corresponding fundamental class $[\overline{\Omega}_I]$. Then in $H^*(\Gr(3,6))$ we have
\begin{equation}\label{gr36ord}
{\tiny\Yvcentermath1\yng(2,1) }
\cdot 
{\tiny\Yvcentermath1\yng(2,1) }
=
{\tiny\Yvcentermath1\yng(3,3)} +
2\ {\tiny\Yvcentermath1\yng(3,2,1)} +
{\tiny\Yvcentermath1\yng(2,2,2)},
\end{equation}
that is, e.g. $\LR_{\{2,4,6\},\{2,4,6\}}^{\{1,3,5\}}=2$.

\smallskip

The natural action of the torus $\T=(\C^*)^n$ on $\C^n$ induces an action of $\T$ on $\Gr(m,n)$. The Schubert cells are invariant, and hence their closures carry a fundamental class in $\T$ equivariant cohomology as well. These classes form a basis of $H^*_{\T}(\Gr(m,n))$ over the ring $H^*_{\T}(\pt)=\Z[z_1,\ldots,z_n]$, where $z_i$ is the first Chern class of the tautological line bundle corresponding to the $i$'th factor of $\T$. Hence the structure constants are polynomials in $z_i$'s. For example, the $\T$ equivariant version of \eqref{gr36ord} now reads
\begin{equation} \label{gr36z}
    \begin{split}
    \MoveEqLeft
\tiny\Yvcentermath1\yng(2,1) \cdot \tiny\Yvcentermath1\yng(2,1) 
 =
     \tiny\Yvcentermath1\yng(3,3) +
2\ \tiny\Yvcentermath1\yng(3,2,1) +
\tiny\Yvcentermath1\yng(2,2,2) +    (2z_5-z_1-z_2) \ \tiny\Yvcentermath1\yng(3,2) 
     \\
    &+
(z_3+z_5+z_6-z_1-z_2-z_4)\ \tiny\Yvcentermath1\yng(3,1,1) +
(z_5-z_6-2z_2)\ \tiny\Yvcentermath1\yng(2,2,1) 
    \\
     & + (z_5-z_4)(z_3+z_5-z_1-z_2) \ \tiny\Yvcentermath1\yng(3,1) +
(z_3-z_2)(z_5+z_6-z_2-z_4) \ \tiny\Yvcentermath1\yng(2,1,1) 
\\
& +(z_5-z_2)^2\tiny\Yvcentermath1\yng(2,2) + (z_5-z_4)(z_5-z_2)(z_3-z_2)\ \tiny\Yvcentermath1\yng(2,1).
  \end{split}
\end{equation}
To recover the non-equivariant version \eqref{gr36ord} from the equivariant \eqref{gr36z}, one needs to substitute all $z_i=0$. 

\bigskip

There are multiple ways of calculating the structure constants presented above, see e.g. \cite{KT} for an effective algorithm tailored to this situation. Now we show an approach which generalizes to the more general settings in Figure \ref{eq:kocka}. This method has two ingredients:

\begin{description} 
\item[(i) Formulas representing fundamental classes] \label{sec:formulas} Consider the $\tb=t_1,\ldots,t_k$ and $\zb=z_1,\ldots,z_n$ variables. Define the rational functions
\[
U_I=\prod_{a=1}^m\prod_{b=i_a+1}^n (z_b-t_a)\prod_{1\leq a<b\leq m} \frac{1}{t_b-t_a}, \qquad
W_I=\Sym_{t_1,\ldots,t_m} (U_I),
\]
and, for a permutation $\sigma\in S_n$ define
\[
W_{\sigma,I}=W_{\sigma^{-1}(I)}(t_1,\ldots,t_m,z_{\sigma(1)},\ldots,z_{\sigma(n)}).
\]
Interpreting 
\begin{itemize}
\item $\tb$ as the Chern roots of the tautological rank $m$ subbundle over $\Gr(m,n)$, and
\item $\zb$ as the tautological Chern roots of the torus $\T$ (cf. Section~\ref{sec:SchLR}), 
\end{itemize}
the function $W_{\id,I}$ represents the fundamental class $[\overline{\Omega}_I]$. 
\item[(ii) Orthogonality] \label{sec:orthogonality}
Define 
\[
\langle f(\tb,\zb),g(\tb,\zb) \rangle = \sum_{K\subset [n]} \frac{f(\zb_K,\zb)g(\zb_K,\zb)}{R_K}, \qquad R_K=\prod_{i\in K} \prod_{j\in [n]-K} (z_j-z_i)
\]
where $|K|=m$, $\zb_K$ is the collection of $z$ variables with index from $K$.
Let $s_0$ be the longest permutation of $n$. Then
\[
\langle W_{\id,I}, W_{s_0,J} \rangle = \delta_{I,J}.
\] 
\end{description}
A direct consequence of the statements in (i) and (ii) is an explicit expression for the structure constants of  $H^*_{\T}(\Gr(m,n))$ with respect to the fundamental classes of Schubert varieties.
\begin{cor} \label{cor:LRord} We have
\begin{equation}\label{eq:LRord}
\LR_{I,J}^K=\langle W_{\id,I} W_{\id,J} , W_{s_0,K} \rangle.
\end{equation}
\end{cor}
The explicit expression in Corollary \ref{cor:LRord} can be coded to a computer, and it can produce expressions like the ones presented in \eqref{gr36z}. It has, however, disadvantages. One of them is the denominators: due to the nature $\langle\ , \rangle$ is defined we obtain the structure constants as a (large) sum of rational functions. Part of the claim is that this rational function in fact simplifies to a polynomial. Yet, such a simplification is usually rather time- and memory-consuming for computers. Even if we are only interested in the non-equivariant structure constants, i.e. the substitution $z_i=0$ in \eqref{eq:LRord}, we must carry out the simplification from rational function to polynomial first, because the denominators of the rational functions are products of $(z_i-z_j)$ factors.  Another disadvantage of the formula for $\LR_{I,J}^{K}$ in Corollary \ref{cor:LRord} is that it does not display known positivity properties of the structure constants. 

\section{The main theorem}\label{sec:main}

The advantage of Corollary \ref{cor:LRord} is that it generalizes to the other vertices in Figure~\ref{eq:kocka}. This feature is the content of the recent development of $\h$-deformed Schubert calculus in cohomology, K theory, and elliptic cohomology.

\begin{thm}[Main Theorem]
Let $\Fl$ be a partial flag variety of type A.
\begin{description}
\item[Formulas for Schubert classes]
There are explicit formulas $W^{\HH}_{\id,I}$,  $W^{\K}_{\id,I}$,  $W^{\E}_{\id,I}$ for the $\h$-deformed Schubert classes in
\begin{equation*}
H_{\T}^*(\Fl), 
\qquad
K_{\T}(\Fl),
\qquad
\Ell_{\T}(\Fl).
\end{equation*}
These are functions in terms of equivariant variables $z_i$, Chern roots of tautological bundles over $\Fl$ called $t^{(i)}_j$, and $\h$ (as well as other parameters in case of $\Ell$).
\item[Orthogonality] 
The given formulas satisfy orthogonality relations for appropriate inner products $\langle \ ,\ \rangle_{\HH}$, $\langle \ ,\ \rangle_{\K}$, $\langle \ ,\ \rangle_{\E}$.
\item[Structure constant formulas] Hence, we have the explicit formulas 
\begin{align*}
\LR_{I,J}^{K}=& \langle W^{\HH}_{\id,I} W^{\HH}_{\id,J}, W^{\HH}_{s_0,K} \rangle_{\HH}, \\
\LR_{I,J}^{K}=& \langle W^{\K}_{\id,I} W^{\K}_{\id,J}, (-\h)^{-\dim_K}\iota[W^{\K}_{s_0,K}] \rangle_{\K}, \\
\LR_{I,J}^{K}= & \langle W^{\E}_{\id,I}W^{\E}_{\id,J}, (\vartheta(\h)/\vartheta'(1))^{\dim_\lambda}  \tau[W^{\E}_{s_0,K}] \rangle_{\E}.
\end{align*}
for the $\h$-deformed Schubert structure constants in $H_{\T}^*$, $K_{\T}$, and $\Ell_{\T}$.
\end{description}
\end{thm}

\noindent The statement of this theorem is deliberately vague, as the details of the theorem are rather technical. The rigorous mathematical meaning of this theorem follows from the explanation of all of its terms throughout the rest of the paper. 

Notations about the partial flag variety $\Fl$ and various structures on it (such as bundles, torus action, Schubert cells and varieties) are set up in Section \ref{sec:Fl}.

After introducing elliptic functions and their trisecant identity in Section \ref{sec:theta}, we present a down-to-earth introduction to the equivariant elliptic cohomology of flag varieties in Section \ref{sec:ellipticH}.

The $\h$-deformed Schubert classes---namely the Chern-Schwartz-MacPher\-son class, the motivic Chern class, and the elliptic class, in $H_{\T}^*, K_{\T}$, and $\Ell_{\T}$---are introduced in Section \ref{sec:classes}.

The formulas  $W^{\HH}_{\id,I}$,  $W^{\K}_{\id,I}$,  $W^{\E}_{\id,I}$ as well as their orthogonality relations are given in Section \ref{sec:weight}.

Some examples for structure constant obtained from the Main Theorem are shown in Section \ref{sec:last}.

\section{The partial flag variety}\label{sec:Fl}

Let $N$ be a positive integer, $\lambda=(\lambda_1,\ldots,\lambda_N)\in \N^N$, and define
\[
\lambda^{(j)}\defeq\sum_{i=1}^j \lambda_i, \qquad n\defeq\lambda^{(N)}=\sum_{i=1}^N\lambda_i.
\]
The partial flag variety $\Fl$ parametrizes nested subspaces
\begin{equation}\label{eq:nest}
\{0\}=V_0 \subset V_1 \subset \ldots \subset V_{N-1}\subset V_N=\C^n
\end{equation}
with $\dim V_j=\lambda^{(i)}$. It is a smooth variety of dimension $\dim_\lambda\defeq\sum_{1\leq i<j \leq N} \lambda_i\lambda_j$. Let us recall the usual structures on $\Fl$.\begin{description}
\item[Bundles] The tautological rank $\lambda^{(i)}$ bundle, whose fiber over the point \eqref{eq:nest} is $V_i$ will be called $\mathcal V_i$. 
\item[Torus action] The standard action of the torus $\T\defeq (\C^*)^n$ on $\C^n$ induces its action on $\Fl$.
\item[Combinatorial gadgets] Consider tuples $I=(I_1,\ldots,I_N)$ where $I_j\subset [n]$, satisfying $|I_j|=\lambda_j$, $I_i\cap I_j=\emptyset$. Their collection will be denoted by $\Il$. 
For example ${\mathcal I}_{(1,2)}=\{ (\{1\},\{2,3\}),$  $(\{2\},\{1,3\}),$ $(\{3\},\{1,2\}) \}$. 
For $I\in \Il$ we will use the notation $I^{(k)}=\bigcup_{s=1}^k I_s = \{i^{(k)}_1<i^{(k)}_2<\ldots<i^{(k)}_{\lambda^{(k)}}\}$.
\item[Torus fixed points] The fixed points $x_I$ of the $\T$ action on $\Fl$ are parametrized by $\Il$: 
\[
x_I=\left( \spa\{\epsilon_i\}_{i\in I_1} \subset \spa\{\epsilon_i\}_{i\in I_1\cup I_2} \subset \ldots \right) \in \Fl,
\]
where $\epsilon_1,\ldots,\epsilon_n$ is the standard basis of $\C^n$.  
\item[Schubert cells] Define the Schubert cell
\[
\Omega_I=\{(V_\bullet)\in \Fl: \dim(V_p \cap \C^q)=|\{i\in I_1\cup\ldots\cup I_p: i\leq q\}|\ \forall p,q\},
\]
where $\C^k=\spa\{\epsilon_1,\ldots,\epsilon_k\}$. We have $x_I\in \Omega_I$ and $\Omega_I$ has dimension 
\[\dim_I\defeq | \cup_{j < k} \{ (a,b)\in I_j\times I_k: a>b\}|.\]
\end{description}

\section{Elliptic functions}\label{sec:theta}

\subsection{Theta functions}
We will use the following version of theta-functions:
\[
\vartheta(x)=(x^{1/2}-x^{-1/2})\prod_{s=1}^\infty(1-q^{s}x)(1-q^{s}/x).
\]
We treat $q\in \C$, $|q|<1$ as a fixed parameter, and will not indicate dependence on it. The function $\vartheta$ is defined on a double cover of $\C$. Theta functions will often appear through
\[
\delta(x,y)\defeq \frac{\vartheta(xy)\vartheta'(1)}{\vartheta(x)\vartheta(y)},
\]
which is meromorphic on $\C^*\times \C^*$.

\smallskip

\begin{rem}
The $q\to 0$ limit we call {\em trigonometric limit} because at $q=0$ the function $\vartheta(x)$ is $\sin(y)$ (up to a constant) in the new variable $x^{1/2}=e^{iy}$. By disregarding the constant factor and denoting the new variable by the same letter as the old one, we say $\vartheta(x)\to \sin(x)$ is our trigonometric limit. The further approximation of $\sin(x)$ with $x$ will be called the {\em rational limit}. The three levels $\vartheta(x) \to \sin(x) \to x$ correspond to the Euler class formulas of line bundles in the three cohomology theories $\Ell$, $K$, $H^*$. \footnote{It is more customary to regard $1-x$ as the K theoretic Euler class of a line bundle, but again, up to a unit ($x$ is invertible in K theory!) this is $\sin(y)$ in a new variable.} Equivalently, the formal group laws of the three theories are (up to constants and change of variables)
\[
(x,y)\mapsto x+y,\qquad
(\sin(x),\sin(y))\mapsto \sin(x+y),\qquad
(\vartheta(x),\vartheta(y))\mapsto \vartheta(xy).
\]
In the three versions the $\delta$-functions (up to constant) are
\[
\frac{x+y}{xy}=\frac{1}{x}+\frac{1}{y}, \qquad
\frac{\sin(x+y)}{\sin(x)\sin(y)}=\cot(x)+\cot(y),\qquad
\frac{\vartheta(xy)}{\vartheta(x)\vartheta(y)}.
\]
Observe that the $\delta$-function separates to a sum of two terms, one depending on $x$ the other on $y$, in the trigonometric and rational limits, but not for theta-functions.
\end{rem}

\subsection{Fay's trisecant identity}
\begin{prop}\cite{F} \cite[Thm. 7.3]{FelRV1}
For variables satisfying $x_1x_2x_3=1$ and $y_1y_2y_3=1$ we have
\begin{equation}\label{eq:Fay}
\delta(x_1,y_2)\delta(x_2,1/y_1)+\delta(x_2,y_3)\delta(x_3,1/y_2)+\delta(x_3,y_1)\delta(x_1,1/y_3)=0.
\end{equation}
\end{prop}
Note that in the trigonometric limit, that is, substituting $\delta(x,y)=\sin(x+y)/(\sin(x)\sin(y))$, identity \eqref{eq:Fay} takes the form
\begin{multline}\label{eq:Fay1}
x_1+x_2+x_3=0, y_1+y_2+y_3=0 \Rightarrow \\
\cot(x_1)\cot(x_2)+\cot(x_2)\cot(x_3)+\cot(x_3)\cot(x_1) = \\
\cot(y_1)\cot(y_2)+\cot(y_2)\cot(y_3)+\cot(y_3)\cot(y_1).
\end{multline}
In the rational limit, that is, substituting $\delta(x,y)=(x+y)/(xy)$, identity \eqref{eq:Fay} takes the form
\begin{multline}\label{eq:Fay2}
x_1+x_2+x_3=0, y_1+y_2+y_3=0 \Rightarrow \\
\frac{1}{x_1x_2}+\frac{1}{x_2x_3}+\frac{1}{x_3x_1}=\frac{1}{y_1y_2}+\frac{1}{y_2y_3}+\frac{1}{y_3y_1}.
\end{multline}
However, in these two limits, more is true. Namely, not only the two sides of \eqref{eq:Fay1} are equal to each other, bot {\em both sides} of \eqref{eq:Fay1} are 0. The same holds for \eqref{eq:Fay2}. The reader is invited to verify that the two sides of \eqref{eq:Fay1} vanish, using high school memories about trigonometric identities. In the elliptic version \eqref{eq:Fay} no such ``separation of $x$ and $y$ variables'' holds. 

 It is worth recording \eqref{eq:Fay1} in ``exponential variables'' ($x_1x_2x_3=y_1y_2y_3=1$):
\begin{multline*}
\frac{1+x_1}{1-x_1}\frac{1+x_2}{1-x_2}+ \frac{1+x_2}{1-x_2}\frac{1+x_3}{1-x_3}+\frac{1+x_3}{1-x_3}\frac{1+x_1}{1-x_1} =\\
\frac{1+y_1}{1-y_1}\frac{1+y_2}{1-y_2}+ \frac{1+y_2}{1-y_2}\frac{1+y_3}{1-y_3}+\frac{1+y_3}{1-y_3}\frac{1+y_1}{1-y_1}, \notag
\end{multline*}
which holds because both sides are equal to $-1$. 

There are various other identities involving theta function, e.g. the ones in \cite[Sect.~2.1]{RTV4} or \cite[Sect.~4.1]{MW} are direct generalizations of \eqref{eq:Fay}.

\section{Equivariant elliptic cohomology of $\Fl$}\label{sec:ellipticH}

This section is an informal introduction to equivariant elliptic cohomology of type~A partial flag manifolds (or, of general so-called GKM spaces). Our general references are \cite{GKV}, \cite[Section 2]{AO},  \cite[Section 4]{FelRV2}, \cite[Section 7]{RTV4}.

Before explaining what we mean by equivariant {\em elliptic} cohomology of $\Fl$ let us revisit its equivariant cohomology and K theory. According to equivariant localization, the restriction maps to the (finitely many) fixed points induce injections of algebras:
\begin{align*}
H^*_{\T}(\Fl)\xhookrightarrow{} & \bigoplus_{x\in \Fl\!\!^{\T}} H^*_{\T}(x)=\bigoplus_{x\in \Fl\!\!^{\T}} \Z[z_1,\ldots,z_n], \\
K_{\T}(\Fl)\xhookrightarrow{} & \bigoplus_{x\in \Fl\!\!^{\T}} K_{\T}(x)=\bigoplus_{x\in \Fl\!\!^{\T}} \Z[z_1^{\pm1},\ldots,z_n^{\pm1}].
\end{align*}
Recall that the $\T$ fixed points $x_I$ of $\Fl$ are parametrized by $\Il$; the map is $f\mapsto (f|_{x_I})_{I\in \Il}$.

Moreover, the image of these injections have the following (so-called GKM-) descriptions. The tuple $(f_I)_{I \in \Fl\!\!^{\T}}$ belongs to the image, if and only if, for ``$(i,j)$-neighboring'' fixed points $x_I$ and $x_J$ the difference of components $f_{I}-f_{J}$ is divisible by $z_i-z_j$. Here ``$(i,j)$-neighboring'' means that $J$ is obtained from $I$ by replacing the numbers $i$ and $j$. Divisibility is meant in the ring of polynomials and in the ring of Laurent polynomials, respectively. It is convenient to rephrase this divisibility condition to
\begin{equation}\label{eq:GKM}
f_I|_{z_i=z_j}=f_J|_{z_i=z_j}\qquad \text{for $(i,j)$-neighboring fixed points $I$ and $J$}.
\end{equation}
Further encoding our descriptions we can say that 
\begin{align*}
H^*_{\T}(\Fl)\xhookrightarrow{} & \bigoplus_{I\in \Il} \text{``natural functions'' on $\C^n$}, \\
K_{\T}(\Fl)\xhookrightarrow{} &\bigoplus_{I\in \Il} \text{``natural functions'' on $(\C^*)^n$},
\end{align*}
such that the image is characterized by \eqref{eq:GKM}. Polynomials and Laurent polynomials are indeed the ``natural functions'' on $\C^n$ and $(\C^*)^n$. This description has a built-in flexibility needed in several applications: in certain studies of $\Fl$ one replaces the coefficient ring $\Z$ with other rings---and this can be achieved by just redefining what ``natural functions'' mean. In some other studies we want to permit some denominators, then again we just need to redefine the notion of ``natural functions'', and still the same description holds. 

In this contexts $\T$ equivariant elliptic cohomology of $\Fl$ would be natural to define by blindly replacing $\C$, or $\C^*$ above with the third 1-dimensional algebraic group, the torus $E=\C^*/(q^{\Z})$ for a fixed $|q|<1$. However, $E$ being compact, there are no functions on $E$ or $E^n$. But there are sections of line bundles. Hence, $\T$ equivariant elliptic cohomology of $\Fl$ is defined by 
\begin{align*}
\Ell_{\T}(\Fl)\xhookrightarrow{} & \bigoplus_{I\in \Il} \text{sections of line bundles on $E^n$},
\end{align*}
with image characterized by \eqref{eq:GKM}. Like above, we use $z_i$ as coordinates, now on $E^n$. Of course one may ask which line bundles are we allowing, and sections of what properties (meromorphic, etc). In our point of view, the answer to those questions reflect different flavors of equivariant elliptic cohomology, similarly to how the exact choice of ``natural functions'' on $\C^n$ or $(\C^*)^n$ resulted different flavors of cohomology and K theory.

\begin{ex} \rm \label{ex:P1}
Let $\Fl=\PPP^1$ and let us permit extra parameters $\h$ and $\mu_1,\mu_2$. Then the ordered pairs
\begin{equation}\label{eq:WP1}
\left( \vartheta(z_2/z_1), 0 \right),\qquad
\left(
\vartheta'(1) \frac{ \vartheta\left(z_2\mu_2/(z_1\mu_1)\right) }{\vartheta(\mu_2/\mu_1)},
\vartheta'(1) \frac{ \vartheta\left(z_1\h/z_2\right)}{\vartheta(\h)}
\right)
\end{equation}
are both equivariant elliptic cohomology classes on $\PPP^1$ (verify the property \eqref{eq:GKM} for both). In fact these two classes will be the $\h$-deformed Schubert classes associated with the Schubert cells $\{(1:0)\}$ and $\PPP^1-\{(1:0)\}$. 
\end{ex}

\begin{rem} \label{rem61}
An easy way to guarantee for a tuple to satisfy condition \eqref{eq:GKM} is to describe the components of the tuple as suitable substitutions of the same function depending on suitable new variables. For example, consider the functions 
\[
\frac{\vartheta(z_1\h\mu_2/(t\mu_1)) \vartheta(z_2/t)}{\vartheta(\h \mu_2/\mu_1)}, \qquad
\vartheta'(1)\frac{ \vartheta(z_1\h/t)\vartheta(z_2\mu_2/(t\mu_1))}{\vartheta(\h)\vartheta(\mu_2/\mu_1)},
\]
and for each one consider the ordered pair of its $t=z_1$ and $t=z_2$ substitutions. We obtain exactly the tuples in the Example above. The very fact that they are $t=z_1$ and $t=z_2$ substitutions of the same function guarantees condition \eqref{eq:GKM}.
\end{rem}

\begin{rem} \rm
In some circumstances specifying the permitted line bundles and the permitted sections is important. For example, if a uniqueness theorem claims that an equivariant elliptic cohomology class is determined by a list of axioms, then one must precisely define which line bundles and what kind of sections are permitted, see \cite[Section 3.5]{AO}, \cite[Appendix A]{FelRV2}, \cite[Section 7.8]{RTV4}. However, if we have some concrete tuple of theta-functions on $E^n$ then we can state that this tuple is an equivariant elliptic cohomology class for the line bundle determined by the transformation properties of the theta functions, as long as the tuple satisfies~\eqref{eq:GKM}.
\end{rem}

\section{$\h$-deformed Schubert classes in $H^*$, $K$, and $\Ell$} \label{sec:classes}

\subsection{$\h$-deformed Schubert class in cohomology: CSM class}

Here we sketch the definition of Chern-Schwartz-MacPherson classes, following \cite{A, AB, O1, O2}, see also \cite{M, AM1, FR}.

Let $F^{\T}(-)$ be the covariant functor of $\T$ invariant constructible functions (on complex algebraic varieties, with an appropriately defined push-forward map using the notion of Euler characteristic). Let $H^{\T}_*(-)$ be the functor of $\T$ equivariant homology as in \cite{edidin}. The $\T$ equivariant MacPherson transformation is the unique natural transformation 
\[
C_*^{\T}: F^{\T}(-) \to H_*^{\T}(-)
\]
satisfying the normalization $C_*^{\T}(\One_X)=c(TX)\cap \mu_X^{\T}$ for smooth projective $X$. Here $\One_X$ is the constant $1$ function on $X$, $c(TX)$ the equivariant total Chern class, and $\mu_X^{\T}$ the equivariant fundamental homology class. If $i:Y\subset X$ is a subvariety of a smooth ambient space $X$ with $\T$ equivariant Poincar\'e duality $\mathcal P$, then we define the (cohomological) Chern-Schwartz-MacPherson (CSM) class $\csm(Y)=\csm(Y\subset X)\defeq \mathcal P (i_*(C_*^{\T}(\One_Y)))\in H^*_{\T}(X)$. \footnote{Although traditionally the class living in {\em homology} is called Chern-Schwartz-MacPherson class---transforming it to the {\em cohomology} of an ambient smooth space is convenient for the purposes of this paper, just like in \cite{O1,O2,FR}.}

The CSM class of $Y\subset X$ is an inhomogeneous cohomology class in $H^*_{\T}(X)$. Its lowest degree component is the fundamental cohomology class $[\overline{Y}]\in H^*_{\T}(X)$. 

It is customary to homogenize it with an extra variable $\h$, making it of homogeneous degree $\dim X$. This version contains the same information as the original $\h=1$ substitution of it, but in some other setups this $\h$-version is more natural (and has its own name ``characteristic cycle class''). For the purpose of this paper we use the $\h$-homogenized one, that is, from now on $\csm(Y)=\csm(Y\subset X)\in H^*(X)[\h]$. In this version it is the coefficient of the highest power of $\h$ in $\csm(Y)$ which equals $[\overline{Y}]$. This justifies our vocabulary of calling the CSM class the {\em $\h$-deformed Schubert class}. 

Let us comment on how one deals with $\csm$ classes in practice. There are three standard approaches. 

The first approach is based directly on the fact that $C_*$ is a natural transformation of functors, and compares the CSM class of $Y$ with the CSM class of its closure and some geometry of the resolution of the closure---taking into account Euler characteristics of fibers.  Typically we arrive at an inclusion-exclusion (sieve) type formula for $\csm(Y)$. This approach can be modified by arranging the inclusion-exclusion argument `upstairs', in the resolution itself.

The second approach is based on the fact that $\csm$, besides the defining axioms, satisfies another strong rigidity property. Its ``Segre version'' $\ssm(Y)=\csm(Y)/c(TX)$ is consistent with  pull-back: $\ssm(f^{-1}(Y))=f^*\ssm(Y)$, for closed $Y$ and sufficiently transversal $f$ to $Y$.

The third approach is that in certain situations the $\csm$ classes satisfy a collection of interpolation constraints, and those constraints uniquely determine them. This approach was triggered by Maulik-Okounkov's notion of {\em cohomological stable envelope}, and was proved in \cite{RV, FR}.

\subsection{$\h$-deformed Schubert class in K theory: motivic Chern class}
The non-$\h$-deformed (equivariant) Schubert calculus has a large literature, going back to \cite{LS, L}, see references in the more recent \cite{GKr}. In this section we follow \cite{BSY, W2, FRW1,  AMSS2} and sketch the definition of the $\h$-deformed Schubert class in K theory: the equivariant motivic Chern class. 

Let $X$ be a quasi-projective complex algebraic $\T$-variety. Let $G_0^{\T}(Var/X)$ denote the Grothendieck group of equivariant varieties and morphisms over $X$, modulo the scissors relation. There is a unique natural transformation 
\[\mC:G_0^{\T}(Var/X) \to K_{\T}(X)[\h]
\]
satisfying\footnote{It is more customary to denote the auxiliary parameter $\h$ by $y$, in accordance with the fact that the integral of the class $\mC[\id_X]$ is the $\chi_y$ genus of $X$ with $\h=y$. Yet, we keep the $\h$ notation to have consistent notation throughout $H^*, K, \Ell$.}
\begin{itemize}
\item functoriality: $\mC[g\circ f]=f_* \mC[g]$, and
\item normalization: $\mC[\id_X]=\lambda_{\h}(T^*X)\defeq \sum \h^i[\Lambda^iT^*X]$ for smooth $X$.
\end{itemize}
For $i:Y\subset X$ we write $\mC[Y]=\mC[i]$. The class $\mC[Y]\in K_{\T}(X)[\h]$ is called the {\em motivic Chern class} of $Y$ in $X$. 

\begin{rem}\rm
For subvarieties $Y$ with mild (so-called Du Bois) singularities, the $\h=1$ substitution recovers the K theoretic fundamental class of $Y$ in $X$. This justifies the name $\h$-deformed class. For $Y$ with non-Du Bois singularities the notion of {\em K theoretic fundamental class} is in fact ambiguous \cite[Section 5]{RSz}, \cite{Fe}, and one may argue that the ``right'' choice for that notion is $\mC[Y]_{\h=1}$.
\end{rem}

Let us comment on how one deals with $\mC$ classes in practice. Similarly to CSM theory there are three different approaches. 

The first one is based directly on the functoriality property: We find a resolution $\tilde{Y} \to \overline{Y}$, and calculate the $\mC$ class of the composition $\tilde{Y}\to X$ using the normalization property. This will not equal $\mC[Y]$, but the difference is supported on the singular locus of $\tilde{Y}\to \overline{Y}$. To find that difference we resolve the singular locus, then the singular locus of that, etc. Finally we arrive at an inclusion-exclusion (sieve) type formula for $\mC[Y]$. This approach can be modified by arranging the inclusion-exclusion argument upstairs, in the resolution itself.

The second approach is based on the fact that $\mC$, besides the defining axioms, satisfies another strong rigidity property. Its ``Segre version'' $\mS[Y]=\mC[Y]/\lambda_{\h}(T^*X)$ is consistent with  pull-back: $\mS[f^{-1}(Y)]=f^*\mS[Y]$, for closed $Y$ and sufficiently transversal $f$ to $Y$.

The third approach is that in certain situations the $\mC$ classes satisfy a collection of interpolation constraints, and those constraints uniquely determine them. This approach was triggered by Okounkov's notion of {\em K theoretic stable envelope}, and is proved in \cite{FRW1, FRW2}.

\subsection{$\h$-deformed Schubert class in elliptic cohomology: the elliptic class} 

The elliptic characteristic class $E(\Omega_I)$ associated with a Schubert cell (in arbitrary $G/P$ type) was defined in \cite{RW, KRW}. This class necessarily depends on the $\h$-variable, as well as a new set of variables $\mu_i$, which are called K\"ahler-, or dynamical variables. Namely,
\[
E(\Omega_I)\defeq \tilde{\mathcal E\!\ell\!\ell}(\overline{\Omega}_I,D_I)
\]
where $\tilde{\mathcal E\!\ell\!\ell}$ is an equivariant and elliptic version of the Borisov-Libgober class \cite{BL1, BL2, BL3, Wa1, Wa2}, and $D_I$ is an appropriate divisor on $\overline{\Omega}_I$ (for details see \cite{RW, KRW}). By setup the Borisov-Libgober class depends on $\h$, and the divisor $D_I$ depends on a character of the corresponding parabolic subgroup $P$. Some of the properties of $E(\Omega_I)$ include the following.
\begin{itemize}
\item $E(\Omega_I)$ specializes to $\mC(\Omega_I)$ and further to $\csm(\Omega_I)$ in the trigonometric, and rational limit of elliptic cohomology.
\item $E(\Omega_I)$ is computable from a resolution of $\overline{\Omega}_I$ through a process similar to the process computing $\mC(\Omega_I), \csm(\Omega_I)$.
\item $E(\Omega_I)$ satisfies and is determined by a small set of axioms which are essentially of interpolation flavor (cf. the interpolation characterization of $\mC(\Omega_I)$, $\csm(\Omega_I)$).
\item In type A the class $E(\Omega_I)$ coincides with the notion of {\em elliptic stable envelope} of \cite{AO}.
\item The switch ``equivariant parameters $\leftrightarrow$ dynamic parameters'' is an incarnation of $d=3$, $\mathcal N=4$ mirror symmetry \cite{RSVZ1, RSVZ2}.
\end{itemize}

\section{Weight functions and their orthogonality relations}\label{sec:weight}

Weight functions, in three flavors---rational, trigonometric, and elliptic---were defined and studied by Tarasov-Varchenko and others in relation with hypergeometric solutions to KZ equations, \cite{TV, RTV1, RTV2, FelRV2, RTV3, RTV4, K}. Here we define weight functions adjusted to our geometric needs, and in Theorem \ref{thm:weightSchubert} we show that they represent $\h$-deformed Schubert classes.

Let $N\in \N$, and $\lambda\in \N^N$. Define the initial sums $\lambda^{(k)}=\sum_{i=1}^k \lambda_i$, and set $n=\lambda^{(N)}$. We will consider functions $W(\tb,\zb,\h)$ and $W(\tb,\zb,\mub,\h)$ in the variables
\[
\tb=(
t^{(1)}_1,\ldots, t^{(1)}_{\lambda^{(1)}},\quad 
t^{(2)}_1,\ldots, t^{(1)}_{\lambda^{(2)}}, \qquad
\ldots \qquad
t^{(N-1)}_1,\ldots, t^{(N-1)}_{\lambda^{(N-1)}}),
\]
\[ 
\zb=(z_1,\ldots,z_n),\qquad \mub=(\mu_1,\ldots,\mu_N),\qquad\text{and}\qquad \h.
\]
When $t^{(N)}_a$ appears in the formulas, it is interpreted as $z_a$. For a function in these variables we define 
\[
\Sym_\lambda (f)  = \Sym_{t^{(1)}}\ldots\Sym_{t^{(N-1)}} (f)
\]
where $\Sym_{t^{(k)}}(g)$  denotes the symmetrization in the $t^{(k)}_*$ variables, i.e. 
\[
\sum_{\sigma\in S_{\lambda^{(k)}}} g\left(t^{(k)}_a\mapsto t^{(k)}_{\sigma(a)}\right).\]

\subsection{Rational weight functions} 
\[
\psi^{\HH}_{I,k,a,b}(x)=\begin{cases}
x+\h & \text{ if } i^{(k+1)}_b<i^{(k)}_a \\
\h & \text{ if } i^{(k+1)}_b=i^{(k)}_a \\
x & \text{ if } i^{(k+1)}_b>i^{(k)}_a, \\
\end{cases}
\]
\[
U_I^{\HH}=\prod_{k=1}^{N-1} 
\left(
\prod_{a=1}^{\lambda^{(k)}}\prod_{b=1}^{\lambda^{(k+1)}} \psi^{\HH}_{I,k,a,b}(t^{(k+1)}_b-t^{(k)}_a) \cdot 
\prod_{a<b\leq \lambda^{(k)}} \frac{1}{t^{(k)}_b-t^{(k)}_a} 
\prod_{b\leq a\leq \lambda^{(k)}} \frac{1}{t^{(k)}_b-t^{(k)}_a+\h} 
\right),
\]
\[
W^{\HH}_I=\Sym_\lambda (U_I^{\HH}), \qquad W^{\HH}_{\sigma,I}=W^{\HH}_{\sigma^{-1}(I)}(\tb,z_{\sigma(1)},\ldots,z_{\sigma(n)},\h)
\qquad (\sigma\in S_n).
\]

\subsection{Trigonometric weight functions}
\[
\psi^{\K}_{I,k,a,b}(x)=\begin{cases}
1+\h x & \text{ if } i^{(k+1)}_b<i^{(k)}_a \\
(1+\h)x & \text{ if } i^{(k+1)}_b=i^{(k)}_a \\
1-x & \text{ if } i^{(k+1)}_b>i^{(k)}_a, \\
\end{cases}
\]
\[
U_I^{\K}=\prod_{k=1}^{N-1} 
\left(
\prod_{a=1}^{\lambda^{(k)}}\prod_{b=1}^{\lambda^{(k+1)}} \psi^{\K}_{I,k,a,b}(t^{(k)}_a/t^{(k+1)}_b) \cdot 
\prod_{a<b\leq \lambda^{(k)}} \frac{1}{1-t^{(k)}_a/t^{(k)}_b} 
\prod_{b\leq a\leq \lambda^{(k)}} \frac{1}{1+\h t^{(k)}_a/t^{(k)}_b} 
\right),
\]
\[
W^{\K}_I=\Sym_\lambda (U_I^{\K}), \qquad W^{\K}_{\sigma,I}=W^{\K}_{\sigma^{-1}(I)}(\tb,z_{\sigma(1)},\ldots,z_{\sigma(n)},\h)
\qquad (\sigma\in S_n).
\]

\subsection{Elliptic weight functions}
Define the integer invariants
\begin{itemize}
\item $p(I,j,i)=|I_j \cap \{1,\ldots,i-1\}$;
\item $j(I,k,a)$ is defined by $i^{(k)}_a\in I_{j(I,k,a)}$,
\end{itemize}
and the functions
\[
\psi^{\E}_{I,k,a,b}(x)=\begin{cases}
\vartheta(x\h)/\vartheta(\h) & \text{ if } i^{(k+1)}_b<i^{(k)}_a \\
\frac{ \vartheta(x\frac{\mu_{k+1}}{\mu_{j(I,k,a)}}\h^{1+p(I,j(I,k,a),i^{(k)}_a)-p(I,k+1,i^{(k)}_a)})}{\vartheta(\phantom{x}\frac{\mu_{k+1}}{\mu_{j(I,k,a)}}\h^{1+p(I,j(I,k,a),i^{(k)}_a)-p(I,k+1,i^{(k)}_a)})}
& \text{ if } i^{(k+1)}_b=i^{(k)}_a \\
\vartheta(x) & \text{ if } i^{(k+1)}_b>i^{(k)}_a, \\
\end{cases}
\]
\[
U_I^{\E}=\vartheta'(1)^{\dim_I}\prod_{k=1}^{N-1} 
\left(
\prod_{a=1}^{\lambda^{(k)}}\prod_{b=1}^{\lambda^{(k+1)}} \psi^{\E}_{I,k,a,b}(t^{(k+1)}_b/t^{(k)}_a) \cdot 
\prod_{a<b\leq \lambda^{(k)}} \frac{1}{\vartheta(t^{(k)}_b/t^{(k)}_a)} 
\prod_{b<a\leq \lambda^{(k)}} \frac{\vartheta(\h)}{\vartheta(\h t^{(k)}_b/t^{(k)}_a)} 
\right),
\]
\[
W^{\E}_I=\Sym_\lambda (U_I^{\E}), \qquad W^{\E}_{\sigma,I}=W^{\E}_{\sigma^{-1}(I)}(\tb,z_{\sigma(1)},\ldots,z_{\sigma(n)},\h,\mub)
\qquad (\sigma\in S_n).
\]

\subsection{Orthogonality}
For $I\in \Il$ let
\begin{align*}
R_I^{\HH}=\prod_{k<l} \prod_{a\in I_k} \prod_{b\in I_l} (z_b-z_a), \qquad\ \ \  & 
Q_I^{\HH}=\prod_{k<l} \prod_{a\in I_k} \prod_{b\in I_l} (z_b-z_a+\h), \\
R_I^{\K}=\prod_{k<l} \prod_{a\in I_k} \prod_{b\in I_l} (1-z_a/z_b), \qquad &
Q_I^{\K}=\prod_{k<l} \prod_{a\in I_k} \prod_{b\in I_l} (1+ z_b/(z_a \h) ), \\
R_I^{\E}=\prod_{k<l} \prod_{a\in I_k} \prod_{b\in I_l} \vartheta(z_b/z_a), \qquad\ \  \  &
Q_I^{\E}=\prod_{k<l} \prod_{a\in I_k} \prod_{b\in I_l} \vartheta(\h z_b/z_a).
\end{align*}

Given $\lambda$ and $I\in \Il$, for a function $f(\tb)$ and $I\in \Il$ let $f(\zb_I)$ denote the result of substituting $t^{(k)}_a=z_{i^{(k)}_a}$, for all $k=1,\ldots,N-1, a=1,\ldots,\lambda^{(k)}$. Define the inner products $\langle\ ,\ \rangle_{\HH}$, $\langle\ ,\ \rangle_{\K}$, $\langle\ ,\ \rangle_{\E}$
\[
\langle f(\tb),g(\tb) \rangle =
\sum_{I\in \Il} \frac{f(\zb_{I})g(\zb_I)}{R_I Q_I}
\]
by using the relevant versions of $R_I$ and $Q_I$ in the denominator. In practice the functions $f, g$ will also depend on other variabes $\zb,\h$ (and $\mub$ in case of $\E$), but the substitution does not affect those.

\begin{thm}[Rational orthogonality]
Let $s_0$ be the longest permutation in $S_n$. We have
\[
\langle W^{\HH}_{\id,I}, W^{\HH}_{s_0,J} \rangle_{\HH} = \delta_{I,J}.
\]
\end{thm}

\begin{thm}[Trigonometric orthogonality]
 Let $\iota[f(\zb,\h)]$ be obtained from the function $f(\zb,\h)$ by substituting $1/t^{(k)}_a$ for $t^{(k)}_a$,  $1/z_i$ for $z_i$ (for all possible indexes) and $1/\h$ for $\h$. We have
\[
\langle W^{\K}_{\id,I}, (-\h)^{-\dim_J}\iota[W^{\K}_{s_0,J}] \rangle_{\K} = \delta_{I,J}.
\]
\end{thm}

\begin{thm}[Elliptic orthogonality]
Let $\tau[f(\zb,\h,\mub)]$ be obtained from the function $f(\zb,\h,\mub)$ by substituting $\h^{\lambda_i}/\mu_i$ for $\mu_i$ (for all $i$). We have
\[
\langle W^{\E}_{\id,I}, (\vartheta(\h)/\vartheta'(1))^{\dim_\lambda}  \tau[W^{\E}_{s_0,J}] \rangle_{\E} = \delta_{I,J}.
\]
\end{thm}

To illustrate the non-triviality of the elliptic orthogonality relations, let us mention that the special case of elliptic orthogonality 
\[
\langle W^{\E}_{\id, (\{3\},\{1,2\})} , (\vartheta(\h)/\vartheta'(1))^2 \tau[ W^{\E}_{s_0,(\{1\},\{2,3\})} ] \rangle_{\E}=0
\]
is equivalent to the trisecant identity \eqref{eq:Fay} with the variables

\centerline{
\begin{tabular}{llll}
$x_1=z_2/z_1$ & & & $y_1=\mu_2/(\mu_1\h)$ \\
$x_2=z_1/z_3$ & & & $y_2=\h$ \\
$x_3=z_3/z_2$ & & & $y_3=\mu_1/\mu_2$. 
\end{tabular}
}

\subsection{Weight functions represent $\h$-deformed Schubert classes}

\begin{thm} \label{thm:weightSchubert}
Interpreting the $t^{(k)}_j$ variables as Chern roots of the tautological bundles $\mathcal V_k$ of rank $\lambda^{(k)}$ and the $z_i$ variables as equivariant variables, the weight functions express the $\h$-deformed Schubert classes:
\[
\begin{tabular}{rcl}
$\csm(\Omega_I)= W^{\HH}_{\id,I}$ & & \cite{RV, FR}, \\
$\mC(\Omega_I)= W^{\K}_{\id,I}$ & & \cite{FRW1}, \\
$E(\Omega_I)= W^{\E}_{\id,I}$ & & \cite{RW, KRW}. 
\end{tabular}
\]
The $\mu_i$ variables in $W^{\E}_{\id,I}$ express the dependence of $E(\Omega_I)$ on the character of the parabolic subgroup. 
\end{thm}

\begin{rem} \label{rem:prob}
The elliptic weight function formulas have singularities at $\h=1$, see the general formulas above, or Remark \ref{rem61} and the Example above it. As mentioned, this is another incarnation of the fact that defining non-$\h$-deformed elliptic classes of Schubert cells in $\Ell(\Fl)$, $\Ell_{\T}(\Fl)$ is problematic.
\end{rem}

\section{Sample Schubert structure constants}\label{sec:last}

Having our Main Theorem (with all its ingredients), it is now only a question of computer power to find Schubert structure constants at each vertex of Figure \ref{eq:kocka}. \footnote{In fact, just 10 of the 12 vertices, see Remark \ref{rem:prob}} In the three subsections below we show some sample calculations in 
\[
H^*(T^*\!\Gr(3,6)) \text{ and } H_{\T}^*(T^*\!\Gr(3,6)), \qquad \Ell_{\T}(T^*\!\PPP^n),\qquad \Ell(T^*\!\PPP^1),
\]
then, in Section \ref{sec:positivity} we discuss questions about these structure constants .

\subsection{Cohomology}

In $H^*(\TGr(3,6))$
\begin{equation} \label{gr36Tstar}
    \begin{split}
    \MoveEqLeft
\tiny\Yvcentermath1\yng(2,1) \cdot \tiny\Yvcentermath1\yng(2,1) 
 =
 \h^9\left(
  \tiny\Yvcentermath1\yng(3,3)+
  2\tiny\Yvcentermath1\yng(3,2,1)+
  \tiny\Yvcentermath1\yng(2,2,2)+
  11 \tiny\Yvcentermath1\yng(3,3,1)+ 
  11 \tiny\Yvcentermath1\yng(3,2,2)+
  46 \tiny\Yvcentermath1\yng(3,3,2)+
 108 \tiny\Yvcentermath1\yng(3,3,3)
\right).
  \end{split}
\end{equation}

Observe that the two extensions \eqref{gr36z} and \eqref{gr36Tstar} of \eqref{gr36ord} go the opposite directions: in one of them the non-zero coefficients extend to ``smaller'' partitions, in the other one to ``larger'' partitions. A combination of the two extensions is, of course, $H_{\T}^*(T^*\!\Gr(3,6))$. Those formulas involve $z_i$ and $\h$, and tend to get very large, yet, for example for $I=(\{2,4,6\},\{1,3,5\})$ ``$=\tiny\Yvcentermath1\yng(2,1)$'', in $H_{\T}^*(T^*\!\Gr(3,6))$ we have 
\[
\begin{split} \MoveEqLeft
\LR_{I,I}^I=
(z_5-z_4)(z_5-z_2)(z_3-z_2) \times \\
& (z_1-z_2+\h)(z_1-z_4+\h)(z_3-z_4+\h)(z_1-z_6+\h)(z_3-z_6+\h)(z_5-z_6+\h).
\end{split}
\]
This coefficient is 0 after substituting $z_i=0$, so the corresponding term is not visible in \eqref{gr36Tstar}. The coefficient of 
$
\h^{\dim_I+\dim_I-\dim_I}=\h^6
$ 
is $(z_5-z_4)(z_5-z_2)(z_3-z_2)$, which turns up in \eqref{gr36z} as $c_{I,I}^{I}$.

\subsection{Equivariant elliptic cohomology}
Consider the elliptic classes for $\Fl=\PPP^n$. For  $k\in [n]$ denote $I_k:=(\{k\},[n]-\{k\})$, and let $\LR_{k,l}^m \defeq \LR_{I_k,I_l}^{I_m}$.  

\begin{thm} 
Let $k\leq l$.  For $m>k$ we have $\LR_{k,l}^m=0$ and  
\[
\LR_{k,l}^{k}=
\vartheta'(1)^{l-1}
\frac{\vartheta(z_l/z_k  \cdot  \mu_2/\mu_1 \cdot \h^{2-l} )}{\vartheta(\mu_2/\mu_1 \cdot \h^{2-l})}
\prod_{i=1}^{l-1}\frac{\vartheta(z_i\h/z_k)}{\vartheta(\h)}
\prod_{i=l+1}^{n}\vartheta(z_i/z_k).
\]
In particular, 
\[
\LR_{k,k}^k=
\vartheta'(1)^{k-1}
\prod_{i=1}^{k-1}\frac{\vartheta(z_i\h/z_k)}{\vartheta(\h)}
 \prod_{i=k+1}^n \vartheta(z_i/z_k).
\]
\end{thm}

\subsection{Non-equivariant elliptic cohomology}
Plugging in $z_i=1$ for all $i$ in equivariant elliptic cohomology formulas yields non-equivariant elliptic cohomology formulas. The actual analysis of occurring functions is intriguing. For example, in $\Ell(T^*{\PPP}^1)$ we obtain
\[
\LR_{(2,1), (2,1)}^{(1,2)}=
\vartheta'(1)^2 \lim_{z_2/z_1\to 1} 
\left( 
\frac{
\frac{ \vartheta(z_2/z_1\cdot \mu_2/\mu_1)}{\vartheta(\mu_2/\mu_1)} - \frac{ \vartheta(z_2/z_1\cdot \h)}{\vartheta(\h)}  
}
{\vartheta(z_2/z_1)}
\right).
\]
Observe that the numerator vanishes at $z_2/z_1=1$, and the denominator has a simple 0 there. Hence the ratio has a removable singularity at $z_2/z_1=1$, and the limit is the value when that singularity is removed.

The two terms in the numerator of the limit above have different {transformation properties} (a.k.a. factors of automorphy) with respect to $z_2\to z_2q$. Hence those terms are sections of different line bundles; therefore, they should not be added unless we choose our vector space to be the direct sum of the vector spaces of sections of different bundles. This is a questionable property of elliptic structure constants which deserves future study. 

\subsection{Positivity?} \label{sec:positivity}
A fundamental feature of both characteristic classes formulas and structure constant formulas in various situations is positivity. For example, the integer structure coefficients in ordinary cohomology (the classical Littlewood-Richardson coefficients) are known to be non-negative. The $z$-polynomial structure coefficients in equivariant cohomology are known to be polynomials of $z_{\text{large}}-z_{\text{small}}$ with non-negative coefficients (see the example in \eqref{gr36z}). Analogous results hold in K theory and equivariant K theory. In the $\h$-deformed worlds, for CSM classes and for motivic Chern classes, some positivity as well as log-concavity results and conjectures can can be found in \cite{AMSS, FRW2}. 

It is reasonable to expect positivity in the $\h$-deformed equivariant elliptic cohomolo\-gy---generalizing the analogous properties in K theory and cohomology. The challenge is to figure out what positivity actually means for an elliptic function. We plan to study these expected elliptic positivity properties in the future. 

\begin{rem}\rm
Besides positivity, other features of characteristic classes and structure constants in various versions of Schubert calculus include stabilization and saturation properties. While hints of stabilization appear in $\h$-deformed cohomology and K theory, nothing is known in elliptic cohomology so far. 
\end{rem}

\end{document}